# Complexity of Computing Quadratic Nonresidues
## N.A. Carella, June 2005

***Abstract:*** This note provides a method for constructing quadratic nonresidues in finite fields of characteristic $p$. It will be shown that there is an effective deterministic polynomial time algorithm for constructing quadratic nonresidues in finite fields.

## 1 Introduction

A quadratic residue $x$ modulo a prime $p$ is simply a square. Every square integer is a square modulo $p$, but every square modulo $p$ is not a square integer. The quadratic symbol

$$(1) \qquad \left(\frac{a}{p}\right) = \begin{cases} -1 & \text{if } x^2 \equiv a \bmod p \text{ has no solution,} \\ 0 & \text{if } a \equiv 0 \bmod p, \\ 1 & \text{if } x^2 \equiv a \bmod p \text{ has a solution,} \end{cases}$$

(or its realization) is the standard method of identifying square and nonsquare elements. A nonzero element $x < p$ is a quadratic residues if the quadratic symbol has the value 1. Otherwise it has the value $-1$ and the element is called a quadratic nonresidue. The problem of determining the complexity of constructing quadratic nonresidues is an open problem in algorithmic number theory. Its calculation is a step in several algorithms, exempli gratia, square roots computing, see [3], construction of irreducible polynomials, see [19], solvability of integer equations, et cetera.

A quadratic nonresidue modulo $p$ is constructible in nondeterministic polynomial time in $O(\log(p)^3)$ bit operations or better. Simply choose $x \in \mathbf{F}_p$ at random and compute the $x^{(p-1)/2} \bmod p$. Since there are $(p-1)/2$ nonzero quadratic nonresidues in $\mathbf{F}_p$, the expected number of trials is 2. In contrast, there is no deterministic polynomial time algorithm, [8, p. 95], [12, p. 15], [16, p. 74] or similar references. In this note it will be shown that the construction of quadratic nonresidues is a deterministic polynomial time operation in any finite fields.

The existence of a deterministic polynomial time algorithm does not necessarily improve the current computational methods of generating quadratic nonresidues nor supercede the random algorithms, as the one described above, because the random algorithms are simpler, easier to implement, and very efficient.

## 2 Preliminaries

This section provides some background information on a few concepts and results used in later sections. The discussions are of limited scopes, and the reader should consult the literature for more detailed results and finer analysis.

*Period Polynomials.*   Let $r = kn + 1$ be prime, and let $g$ be a primitive root modulo $r$. The *Gauss period* $\eta_j$ and *coperiod* $g_j$ are defined by the incomplete exponential sums

$$(2) \qquad \eta_j = \sum_{x=0}^{k-1} e^{i2\pi g^{nx+j}/r}, \ \text{ and } \ g_j = \sum_{x=0}^{r-1} e^{i2\pi g^{j}x^{n}/r},$$

for $j = 0, 1, \ldots, n - 1$.

Nontrivial periods of order $n > 1$ exist if and only if $r = kn + 1$ is a prime or $kn = \varphi(r)$, where $\varphi(r)$ is the totient function, and $r \in \mathbb{N}$ is an integer.

The *period polynomial* and its associated *coperiod polynomial* are defined by

$$(3) \qquad \psi_r(x) = (x - \eta_0)(x - \eta_1) \cdots (x - \eta_{n-1}) \ \text{ and } \ \theta_r(x) = (x - g_0)(x - g_1) \cdots (x - g_{n-1}),$$

where the degree $\deg(\psi_r) = n$. The period polynomial of maximal degree $kn = \varphi(r)$ is called the $r$th cyclotomic polynomial $\Phi_r(x)$.

From the defining sums for $\eta_j$ and $g_j$, it is clear that $g_0 = 1 + n\eta_0$, $g_1 = 1 + n\eta_1$, ..., $g_{n-1} = 1 + n\eta_{n-1}$. Accordingly, the period and the coperiod polynomials are equivalent up to a linear change of variable. The linking formulae are

$$(4) \qquad \psi_r(x) = \prod_{i=0}^{n-1} (x - \eta_i) = n^{-n} \theta_r(nx+1), \ \text{ and } \ \theta_r(x) = \prod_{i=0}^{n-1} (x - g_i) = n^n \psi_r((x-1)/n).$$

In some cases it is easier to compute $\theta_r(x)$, and then $\psi_r(x)$ using the above formula.

The coefficients of period polynomials of small degrees can be determined in term of the parameters of the quadratic partitions of the integers $r = a^2 + ub^2$, $a, b, u \in \mathbb{Z}$ or similar equations. Some of them are completely determined. The first two even degrees are listed here.

$n = 2$. If $r = 2k + 1$ is a prime, then

$$(5) \qquad \psi_r(x) = x^2 + x + (1 - (-1)^{(r-1)/2}r)/4, \text{ and } \theta_r(x) = x^2 - (-1)^{(r-1)/2}r.$$

$n = 4$. If $r = 4k + 1 = a^2 + b^2$ is a prime, $a \equiv 1 \bmod 4$, and $k$ is even, then

$$(6) \qquad \psi_r(x) = x^4 + x^3 - \frac{3r-3}{8}x^2 - \frac{(2a+3)r-1}{16}x - \frac{(4a^2+8a-r+7)r-1}{256},$$

and $\theta_r(x) = (x^2 - r)^2 - 4r(x - a)^2$. Otherwise

$$(7) \qquad \psi_r(x) = x^4 + x^3 - \frac{r-3}{8}x^2 - \frac{r-(a+1)2r-1}{16}x - \frac{(4a^2-8a-9r-2)r-1}{256},$$

and $\theta_r(x) = (x^2 + 3r)^2 - 4r(x - a)^2$, see [4, p. 326].

*Irreducibility*

The period polynomials $\psi_r(x)$ of degree $\deg(\psi_r) = n$ are irreducible over the rational numbers $\mathbb{Q}$, and factor as $\psi_r(x) = \psi_{r,1}(x)\psi_{r,2}(x) \cdots \psi_{r,s}(x)$ over the finite field $\mathbf{F}_q$. There are several factorization patterns, for example, a factorization pattern can have $s = n/\mathrm{ord}_r(q)$ irreducible factors $\psi_{r,i}(x)$ over $\mathbf{F}_q$ of degree $d_i = \deg(\psi_r) = \mathrm{ord}_r(q)$. The factorization patterns are specified by a well established result.

**Theorem 1.** Let $r$ be a squarefree integer, and $q$ a prime with $\gcd(r, q) = 1$. Let $K$ be a subgroup of index $e$ in $\mathbb{Z}_r$. Let $d$ be the smallest integer for which $q^d \in K$. Then $d = \mathrm{lcm}(d_1, d_2, \ldots, d_s)$, where $d_i$ is the degree of the $i$th irreducible factor $\psi_{r,i}(x)$ of the period polynomial $\psi_r(x) \in \mathbf{F}_q[x]$.

The proof, derived from the Kummer-Dedekind factorization theorem, appears in [10].

Under the right conditions stated in Theorem 1, the irreducible period polynomials $\psi_r(x) \in \mathbb{Q}[x]$ remain irreducible in $\mathbf{F}_q[x]$. A recent application of this property to primality testing appears in [15]. To state the result the following definition is needed. A *period system* for an integer $n > 1$ is a sequence $(r_1, q_1), \ldots, (r_k, q_k)$ of ordered pairs that satisfies the following properties:
(a) $r_1, r_2, \ldots, r_k$, are primes, (b) $q_1, q_2, \ldots, q_k$, are pairwise relatively primes, (c) $q_i > 1$ divides $r_i - 1$, and the order of $n^{(r_i-1)/q_i}$ modulo $r_i$ is $q_i$.

**Theorem 2.** ([15]) There is a deterministic algorithm such that for each integer $m > 0$ the algorithm produces an integer $D_m$ and further, for each integer $n > 1$, and each integer $D$ with $D > D_m$ and $D > (\log n)^{11/6+1/m}$, the algorithm finds a period system $(r_1, q_1), \ldots, (r_k, q_k)$ for $n$ with each $r_i < D^{6/11}$ and each $q_i < D^{3/11}$, with $D \le q_1 \cdots q_k < 4D$, and $k = O((\log\log D)^2)$. The running time of this algorithm is $O(D^{12/11})$. The implied constant may depend on the choice of $m$.

This result is utilized to construct irreducible period polynomials over $\mathbf{F}_q$.

*Complexity.* An irreducible period polynomial $\psi_r(x)$ of degree $d = \deg(\psi_r)$ over $\mathbf{F}_q$ is computable in deterministic polynomial time in $(dr\log q)$ finite field operations. The algorithm appears in [1], [15], and a complete analysis is given in the latter. A simplified version of the algorithm is given here. In this algorithm and in applications of period polynomials to finite fields it is convenient to rewrite the periods as

$$(8) \qquad \eta = \sum_{x \in K} e^{i2\pi x/r}, \; \eta_j = \eta^{q^j},$$

and the period polynomial as

$$(9) \qquad \psi_r(x) = (x - \eta)(x - \eta^q)(x - \eta^{q^2})\cdots(x - \eta^{q^{d-1}}),$$

where $r = de + 1$, $K = \{ x^e : 0 \le x < r/d \}$, and $\text{ord}_r(q) = d$.

***Algorithm* 3.** ([15]) Input: A list of periods $\eta_0, \eta_1, \ldots, \eta_{d-1}$.
Output: An irreducible period polynomial $\psi_r(x)$ of degree $d$.
Step 1. Equal degrees multiply: multiply the linear factors in pairs, the resulting quadratic factors in pairs, etc. For example, $f_0(x) = (x - \eta_0)(x - \eta_1)$, $f_1(x) = (x - \eta_2)(x - \eta_3)$, $f_2(x) = (x - \eta_4)(x - \eta_5)$, …, $f_0(x)f_1(x)$, $f_2(x)f_3(x)$, $f_4(x)f_5(x)$, …, $f_0(x)f_1(x)f_2(x)f_3(x)$, ….
Step 2. Unequal degrees multiply: multiply the remaining unequal degrees terms to complete the product $\psi_r(x) = (x - \eta_0)(x - \eta_1) \cdots (x - \eta_{d-1})$.

***Remark* 4.** It appears that Algorithm 3 can be used to factor the $r$th cyclotomic polynomial $\Phi_r(x)$ over a finite field $\mathbf{F}_q$ in deterministic polynomial time. For example, apply the algorithm to the linear factors $x - e^{i2\pi/r}$, $x - e^{i2\pi q/r}$, $x - e^{i2\pi q^2/r}$, ..., $x - e^{i2\pi q^{d-1}/r}$, for each $a \ne 0$ such that $\gcd(a, r-1) = 1$. And probably not much work is needed to extend it to all polynomials.

*The Trace Matrix and Its Determinant.* A basis $\{ \alpha_0, \alpha_1, \ldots, \alpha_{n-1} \}$ of the finite field $\mathbf{F}_{q^n}$ over $\mathbf{F}_q$ is a subset of elements whose linear span $\{ a_0\alpha_0 + a_1\alpha_1 + \cdots + a_{n-1}\alpha_{n-1} : a_i \in \mathbf{F}_q \}$ represents the finite field $\mathbf{F}_{q^n}$ as a vector space over $\mathbf{F}_q$. The most common basis is the *power basis* $\{ \alpha_0 = 1, \alpha_1 = x, \alpha_2 = x^2, \ldots, \alpha_{n-1} = x^{n-1} \}$. Given a basis, the *trace matrix* $T$ is defined by the $n \times n$ matrix $T = ( Tr(\alpha_i\alpha_j) )$, where $0 \le i, j \le n - 1$, and $Tr : \mathbf{F}_{q^n} \to \mathbf{F}_q$ is the trace function $Tr(x) = x + x^q + x^{q^2} + \cdots + x^{q^{n-1}}$. The matrix $T$ is a circulant matrix or nearly circulant, so one or two of its rows are sufficient to specify it

***Theorem* 5**. The determinant of the matrix $T$ is of the form

$$(10) \qquad \det(T) = \begin{cases} \text{Quadratic residue in } \mathbf{F}_q & \text{if } n \text{ is odd,} \\ \text{Quadratic nonresidue in } \mathbf{F}_q & \text{if } n \text{ is even.} \end{cases}$$

This result is derived from Theorem 1 and Lemma 1 in [14], and the proof is essentially the same.

*Complexity.* The trace $Tr$ is computable in $O(n^3 \log(q))$ basic finite field $\mathbf{F}_q$ operations or better, see the literature. In the case of matrices, the expanded expression of the determinant of an $n \times n$ matrix has $n!$ terms. Nevertheless it can be computed in running time $O(n^4)$ or better, (by reducing the matrix to a triangular matrix and taking the product of the diagonal entries or some other methods), see [5] for details. From these observations it follow that the determinant of the $n \times n$ trace matrix is computable in about $O(n^4 \log(q))$ basic finite field $\mathbf{F}_q$ operations.

## 3 The Current Time Complexity Of Quadratic Nonresidues

The traditional approach to the determination of a quadratic nonresidue attempts to find the least positive quadratic nonresidue $n_p$ modulo $p$. In the 1800's it was determined that $n_p < p^{1/2} + 1$, [11, Art. 129]. The modern methods of investigating the existence of the least quadratic nonresidues are by means of exponential sums and/or $L$-functions analysis. For example, estimating the least integer $N > 1$ for which exponential sum

$$(11) \qquad \sum_{x=1}^{N} \left( \frac{x}{p} \right) < N.$$

Around the 1950's several authors used exponential sums analysis to reduce the estimate $n_p < p^{1/2} + 1$ to $n_p < p^{1/(4\sqrt{e}) + \varepsilon}$, where $e = 2.71...$, and $\varepsilon > 0$, see [7]. The well known Vinogradov's conjecture claims that $n_p = O(p^\varepsilon)$, $\varepsilon > 0$. Further, on the basis of the extended Riemann hypothesis for $L$-functions it has been established that $n_p = O(\log(p)^2)$, see [2], and [6] for its generalization. This last conditional result implies that a quadratic nonresidue is constructible in fewer than $O(\log(p)^5)$ bit operations. In the other direction, it has been unconditionally proven that there are infinitely many primes for which $n_p \geq c_0 \log(p) \log\log\log(p)$, $c_0$ an absolute constant, see [13]. Primes of the form $p = 2^{u_0} \cdot 3^{u_1} \cdot 5^{u_2} \cdots p_n^{u_n} + 1$, where $p_n$ is the $n$th prime, (or similar forms), have the largest least quadratic nonresidues.

Using the Quadratic Reciprocity Law

$$(12) \qquad \left( \frac{p}{q} \right) \left( \frac{q}{p} \right) = (-1)^{(p-1)(q-1)/4}, \qquad \left( \frac{2}{p} \right) = (-1)^{(p^2-1)/8},$$

it is straightforward to obtain the following quadratic nonresidues:

(i) If $p = 8n + 3$, then $q = -1$ and $2$ are quadratic nonresidues.
(ii) If $p = 8n + 5$, then $q = 2$ is a quadratic nonresidue.
(iii) If $p = 8n + 7$, then $q = -1$ is a quadratic nonresidue.

These are quadratic nonresidues $z$ of least or the least absolute values $|z| \leq 2$. Thus the primes in the arithmetic progression $p = 8n + 1$ are the only remaining ones that do not have known quadratic nonresidues. The quadratic reciprocity law can be used to produces quadratic nonresidues of any equivalence class of primes $p = 2^k n + 1$, $k > 2$, $n$ odd. However, this is not practical: repeatedly using (12) to obtain quadratic nonresidues for $k = 3, 4, 5, \ldots$ requires infinitely many equivalence classes to cover all the primes.

Thus the best unconditional deterministic algorithm has exponential time complexity, and as far as I am aware there is no subexponential algorithm either.

Note that a single quadratic nonresidue $z$ is sufficient to generate the entire set of quadratic nonresidues. This is accomplished by multiplication by squares: $s \rightarrow sz$, where $s$ runs over the set of quadratic residues $Q = \{ x^2 \bmod p : 0 < x < p/2 \}$. Moreover, the entire sequences of consecutive pairs and equally spaced pairs are generated by

$$(11) \qquad z_s = \frac{(sz \pm 1)^2}{4sz}, \ \ z_s \mp 1, \ \ \text{and} \ \ z_s = \frac{(sz \pm v)^2}{4sz}, \ \ z_s \mp v.$$

## 4 General Result

This section presents a general method of generating quadratic nonresidues in any finite fields $\mathbf{F}_q$, $q = p^{2t+1}$, $t \geq 0$. In the case of prime finite fields $\mathbf{F}_p$, the ranks or absolute values of the quadratic nonresidues $z \in \mathbf{F}_q$ are uncertain. But it is possible that $|z| \leq r$ since period polynomials of order $r \leq O(\log(q)^c)$ are involved.

The existence of small integers $r \leq O(\log(q)^c)$ for which the prime power $q = p^{2t+1}$ have even orders modulo $r$ is a crucial part of the result. In some cases, even order $d = \mathrm{ord}_r(q) = 2l$ implies that $r$ is a quadratic nonresidue modulo $q$, but in general this is not true. For example, any prime powers in the arithmetic progressions $q = rx + a = 5x + 4, 9x + 4, 9x + 7, x \geq 0$, etc.

**Theorem 6**. There is a deterministic polynomial time algorithm that constructs quadratic nonresidues in any finite field $\mathbf{F}_q$, $q = p^{2t+1}$, $t \geq 0$.

Proof: Using a weaker modified version of Theorem 2, it follows that there exists an integer $r < 4\log(q)^2$ of multiplicative order $d = \mathrm{ord}_r(q) = 2l \geq 2$. And an irreducible period polynomial $\psi_r(x) \in \mathbf{F}_q[x]$ of degree $d = \deg(\psi_r(x))$ is generated by Algorithm 3. Therefore by Theorem 4, the determinant $\det(T) \in \mathbb{Z}$ of the trace matrix $T = ( Tr(x^{i+j}) )$ is a quadratic nonresidue in $\mathbf{F}_q$.

The time complexity of the algorithm is determined by the running time complexities of computing the period polynomial and computing the determinant of the trace matrix. Both of these have deterministic polynomial time complexities, see section 2. ∎

***Remark* 7**. The Lenstra and Pomerance Theorem is a stronger result than required in Theorem 6. In fact a single integers pair ($r = r_i$, $d = q_i$) such that the order $d = \mathrm{ord}_r(q) = 2l \geq 2$ would do, and this is a weaker requirement since the density of integers $r \leq x$ such that $\mathrm{ord}_r(q) = 2l \geq 2$ is better than $1/2$, see [20].

***Example* 8**. For any prime power $q = p^{2t+1} = 3k + 2$, the polynomial $f(x) = x^2 + x + 1$ is irreducible over $\mathbf{F}_q$, and the set $\{\, 1, x \,\}$ is a basis of $\mathbf{F}_{q^2}$ over $\mathbf{F}_q$. The trace matrix is

$$(12) \qquad T = \begin{bmatrix} Tr(1) & Tr(x) \\ Tr(x) & Tr(x^2) \end{bmatrix} = \begin{bmatrix} 2 & -1 \\ -1 & -1 \end{bmatrix}.$$

Thus its determinant $det(T) = -3$, which coincides with the discriminant of $f(x)$, is a quadratic nonresidue in $\mathbf{F}_q$. Likewise, for $q = p^{2t+1} = 4k + 3$, the polynomial $f(x) = x^2 + 1$ is irreducible over $\mathbf{F}_q$, and $det(T) = -1$ is a quadratic nonresidue in $\mathbf{F}_q$. These solve the quadratic nonresidue problem for all the prime powers $q = p^{2t+1} \neq 12k + 1$. Observe that the same calculations for the prime powers $q = 12k + 1$ with a single irreducible polynomial $x^2 + ax + b$ seems to be impossible.

The general result works with any class of irreducible polynomials of even degrees. Moreover, exploiting the structure of a specific class of polynomials can transform the basic algorithm into a distinct algorithm. In the case of period polynomials, the technique of projecting a higher degree period down to a lower degree period achieves substantial simplification and transparency in the proof over the general algorithm. This is demonstrated below.

## 5 Special Result

A simpler and more practical algorithm for generating quadratic nonresidues in almost any finite fields $\mathbf{F}_q$, $q = p^{2t+1}$, $t \geq 0$, is considered here. The technique is quite general and applies to all period polynomials. Most importantly, it does not rely on any overly complicated theory, and the algorithm has low complexity.

***Theorem* 9**. Let $q \neq 2^{u_0} \cdot 3^{u_1} \cdot 5^{u_2} \cdot 17^{u_3} \cdot 257^{u_4} \cdot 65537^{u_5} \cdot N + 1$, and $u_0 > c_1 \log\log(q)$, $c_1 > 0$ constant, where $\gcd(2 \cdot 3 \cdot 5 \cdot 17 \cdot 257 \cdot 65537, N) = 1$. Then a quadratic nonresidue in $\mathbf{F}_q$, $q = p^{2t+1}$, $t \geq 0$, can be generated in deterministic polynomial time in about $O((\log q)^2)$ finite field operations.

Proof: Let $k = c_1 \log\log(q)$, and define the subset of integers

$$(15) \qquad R = \{\, 2^{v_0} \cdot 3^{v_1} \cdot 5^{v_2} \cdot 17^{v_3} \cdot 257^{v_4} \cdot 65537^{v_5} : 0 \leq v_0 \leq k, 0 \leq v_i \leq 1 \,\}.$$

By hypothesis $q - 1 \not\equiv 0 \bmod r$ for some integer $r \in R$. Under this condition the order of $q$ modulo $r$ is even $\mathrm{ord}_r(q) = 2^m$, $1 \leq m \leq k + 16$. This in turns implies that there is an $r$th root of unity of $\omega \in \mathbf{F}_{q^{2^m}}$ such that $\omega \notin \mathbf{F}_q$. In addition, the projection of $\omega$ into a quadratic extension $\mathbf{F}_{q^2}$ of $\mathbf{F}_q$ is precisely

(14) $$\tau = Tr(\omega) = \omega + \omega^{q^2} + \omega^{q^{2 \cdot 2}} + \cdots + \omega^{q^{2(2^{m-1}-1)}}.$$

Clearly $\tau^q \neq \tau$, but $\tau^{q^2} = \tau$. Thus the period polynomial $\psi_r(x) = x^2 - (\tau + \tau^q)x + \tau\tau^q \in \mathbb{Z}[x]$ remains irreducible in $\mathbf{F}_q[x]$, and its discriminant $D = \tau^2 + \tau^{2q} - 2\tau\tau^q \in \mathbb{Z}$ is a quadratic nonresidue in $\mathbf{F}_q$. Note that $\psi_r(x) \in \mathbb{Z}[x]$ is computed using $\omega = e^{i2\pi/r}$ in Algorithm 3, ∎

The small subset of integers $R$ solves the quadratic nonresidues problem for almost every finite field $\mathbf{F}_q$, $q = p^{2t+1}$, $t \geq 0$, this includes large proportions of the difficult cases

(i) Prime powers $q = p^{2t+1} = 2^k n + 1$, with arbitrary high power of 2,
(ii) Prime powers $q = 2^{u_0} \cdot 3^{u_1} \cdot 5^{u_2} \cdots p_n^{u_n} + 1$, where $p_n$ is the $n$th prime.

**Example 10.** Given any prime power $q = p^{2t+1} = 3n + 2$, the order $\mathrm{ord}_3(q) = 2$, since the prime $r = 3$ does not divides $q - 1$. Accordingly, there is a cube root of unity $\omega \in \mathbf{F}_{q^2}$ such that $\omega \notin \mathbf{F}_q$, and $\tau = \omega$.

Thus the period polynomial $\psi_r(x) = x^2 - (\omega + \omega^q)x + \omega\omega^q = x^2 + x + 1$ is irreducible over $\mathbf{F}_q$, and its discriminant $D = \omega^2 + \omega^{2q} - 2\omega\omega^q = -3$ is a quadratic nonresidue in $\mathbf{F}_q$. The algorithm of Theorem 9 performs similar calculations for $r = 3, 4, 5, 6, \dots$.

## 5 Some Applications

The square roots of quadratic elements in finite fields are efficiently determined using nondeterministic polynomial time algorithms, less than $O(\log(p)^4)$ finite fields operations. These algorithms are fast and easy to implement, see [3], [16], [17], etc. for descriptions. On the other hand, the current best deterministic polynomial time algorithm has a high running time complexity and it is complicated. This algorithm runs in $O(|x|^{1/2+\varepsilon}\log(p)^9)$ bit operations, any $\varepsilon > 0$. Moreover, there is a dependence on the absolute value of the argument $x$, but for $p \not\equiv 1 \bmod 16$, it is not dependent on $|x|$, see [18].

*Corollary* 11. There is a deterministic polynomial time algorithm that computes square roots in finite fields $\mathbf{F}_q$, $q = p^{2t+1}$, $t \geq 0$. Specifically, if $p = 2^k n + 1$, $k \leq c \log\log(q)$, $c > 0$ constant, and the equation $x^2 \equiv a \bmod p$ is solvable, then the roots are given by

(16) $$\sqrt{a} = \pm 2^{2-k} a^{(p+2^k-1)2^{1-k}} \sum_{i=0}^{2^{k-1}-1} (z^{(2i+1)n} - 1)^{-1} a^{ni},$$

where $z$ is a quadratic nonresidue. In the other case $k > c \log\log(q)$ it is computed recursively.

The last formula is due to Cipolla, see [9, Vol. I, p. 220].

***Corollary* 13**.  Irreducible polynomials $x^{2^e} + a \in \mathbf{F}_q[x]$, $q = p^{2t+1}$, $t \geq 0$, of degree $2^e \geq 2$ are constructible in deterministic polynomial time.

The specific details of generating a polynomial of degree $2^e$ from an irreducible quadratic polynomial $x^2 + a \in \mathbf{F}_q[x]$ are given in [19, p. 439.]. Moreover, a single polynomial $f(x)$ of degree $n$ serves as a generator of all the other irreducible polynomials of the same degree by means of the linear fractional transformations

(17) $$(cx + d)^n f((ax + b)/(cx + d)) = a_n x^n + a_{n-1} x^{n-1} + \cdots + a_1 x + a_0,$$

where $ad - bc \neq 0$.